\documentclass{article}
\usepackage{amsmath,amssymb}
\usepackage{graphicx}

\begin{document}

\title{A Common Fixed Point Theorem For Weakly Compatible Mappings Satisfying Property (E. A)\thanks{%
Mathematics Subject Classifications: 47H10,
54H25.}}
%\date{{\small 20 July 2000}}
\author{G. V. R. Babu \thanks{%
Department of Mathematics,  Andhra University, Visakhapatnam-530 003,
Andhra Pradesh, India; e-mail address: gvr$_{-}$babu@hotmail.com}\ , G. N. Alemayehu\thanks{%
Department of Mathematics,  Andhra University, Visakhapatnam-530
003, Andhra Pradesh, India; Permanent address: Department of
Mathematics,
 Jimma University, Jimma, P.O.Box 378, Ethiopia; e-mail address:
alemg1972@gmail.com}\ } \maketitle

\begin{abstract}
The aim of this paper is to prove the existence of common fixed points 
for a pair of weakly compatible selfmaps satisfying weakly contractive condition and property (E. A). In this context, first we modify Beg and Abbas theorem (\cite{Beg}, Theorem 2.5). Our results improve the corresponding results of Beg and Abbas \cite{Beg}. 
\end{abstract}

\section{Introduction}

In 2002, Aamri and Moutawakil \cite{Aamri}  introduced the notion of
property (E. A). There are a number of results (Aliouche \cite
{Aliouche}, Imdad \emph{et. al.} \cite{Imdad}, Liu \emph{et. al.}
\cite{Liu}, Pathak \emph{et. al.} \cite{Pathak}) that use this
concept to prove existence results in common fixed point theory.

Thoughout this paper, $(X,d)$ denotes a metric space; and $f$ and
$T$ are selfmaps of $X$.

\smallskip

DEFINITION 1.1. The pair $(f, T)$ is said to  \\(i) be
\emph{compatible} (Jungck \cite{Jungck1}) if $\lim\limits_{n\to
\infty} d(fTx_{n}, Tfx_{n})=0$, whenever $\{x_{n}\}$ is a sequence
$~~~~~$in $X$ such that $\lim\limits_{n\to \infty}fx_{n}=\lim\limits_{n\to
\infty}Tx_{n}=t,$ for some $t$ in $X$. \\ (ii) be
\emph{noncompatible}  if there is at least one sequence $\{x_{n}\}$
in $X$ such that $\lim\limits_{n\to \infty}fx_{n}\\$$~~~~~=\lim\limits_{n\to
\infty}Tx_{n}=t,$ for some $t$ in $X$, but $\lim\limits_{n\to
\infty} d(fTx_{n}, Tfx_{n})$ is either non-zero or \\$~~~~~~$non-existent.
\\(iii) \emph{satisfy property} (E. A) (Aamri and El
Moutawakil \cite{Aamri}) if there exists a sequence $~~~~~~$$\{x_{n}\}$ in
$X$ such that $\lim\limits_{n\to \infty}fx_{n}=\lim\limits_{n\to
\infty}Tx_{n}=t,$ for some $t$ in $X$.\\
 (iv) be \emph{weakly compatible} (Jungck
\cite{Jungck2}) if $Tfx=fTx$ whenever
 $fx=Tx,~x\in X$.

\smallskip
\smallskip

REMARK 1.2. Every pair of noncompatible selfmaps of a metric space
$(X,d)$ satisfies property (E. A), but its converse need not be
true as shown by the following example.

\smallskip
\smallskip

EXAMPLE 1.3. Let $X=[0,1)$ with the usual metric. We define mappings
$f$ and $T$ on $X$ by

\[\begin{array}{lr}
 f(x)=\left\{\begin{array}{ll}
                  ~~\frac{2}{3} & \text{if}~0\leq x < \frac{2}{3} \\\\
                 1-\frac{1}{2}x & \text{if} ~\frac{2}{3}\leq
  x< 1 \end{array}\right.\nonumber &\mbox{}\text{ and } ~~~~  T(x)=\left\{\begin{array}{ll}
                  ~~\frac{2}{3} & \text{if}~0\leq x < \frac{2}{3} \\\\
                 \frac{4}{3}-x & \text{if} ~\frac{2}{3}\leq
  x< 1 .\end{array}\right.\nonumber
\end{array}\]

\smallskip

Then the pair $(f,T)$ is compatible on $X$, for if $\{x_{n}\}$ is a
sequence in [0,1) with $\lim \limits_{n\rightarrow
\infty}fx_{n}=\lim \limits_{n\rightarrow \infty}Tx_{n}=z\in X$, then
$z=\frac{2}{3}$ and $\lim \limits_{n\rightarrow \infty}Tfx_{n}=\lim
\limits_{n\rightarrow \infty}fTx_{n}=\frac{2}{3}$ so that $\lim
\limits_{n\rightarrow \infty}d(Tfx_{n},fTx_{n})=0$. Hence the pair
$(f,T)$ is not noncompatible on $X$. We observe that $(f,T)$
satisfies property $(E,A)$.
\smallskip

REMARK 1.4. (i) Weak compatibility and property (E. A)
are independent to each other(Pathak \emph{et. al.} \cite{Pathak}).\\
(ii) Every compatible pair is weakly compatible but its converse
need not be true \\(Jungck \emph{et. al.} \cite{Jungck2}).

\smallskip

 Throughout this paper, we denote $R_{+}=[0,\infty )$; and  $N$,  the set of all natural numbers, and
$\Phi =\{\varphi|~\varphi:R_{+}\rightarrow R_{+}~ is ~continuous,
~\varphi(0)=0, ~\varphi(t)> 0$ for $t> 0\}$.

In 2007, Beg and Abbas \cite{Beg}  established the following
existence theorem of common fixed points of a pair of selfmaps.

\smallskip

THEOREM 1.5. (Beg and Abbas \cite{Beg}, Theorem 2.5). Let $(X,d)$
be a metric space and let $T,~f:X\rightarrow X$ be weakly compatible
selfmaps. Assume that there exists a monotone increasing $\varphi
\in \Phi$ with $\lim \limits_{t\rightarrow \infty} \varphi(t)=\infty
$ such that for all $x,y\in X$,

\[d(Tx,Ty)\leq d(fx,fy)-\varphi(d(fx,fy)).\] If $T(X)\subseteq f(X)$
and $f(X)$ is a complete subspace of $X$, then $f$ and $T$ have a
unique fixed point in $X$.

\smallskip

\smallskip

DEFINITION 1.6. A selfmap $T:X\rightarrow X$ is said to be
\emph{weakly contractive with respect to a selfmap $f:X\rightarrow
X$} if there exists a $\varphi \in \Phi$ such that for all $x,y\in
X$,
\begin{equation}
d(Tx,Ty)\leq d(fx,fy)-\varphi(d(fx,fy)). \label{1}
\end{equation}

\smallskip

The aim of this paper is to give a modified version of Theorem 1.5
by relaxing the conditions `$\varphi$ is monotonically increasing
and $\lim \limits_{t\rightarrow \infty} \varphi(t)=\infty $'.

Further we prove the existence of common fixed points for a pair of
weakly compatible selfmaps satisfying  weakly contractive condition
and property (E. A).

\section{A Modified Version Of Beg and Abbas Theorem}

The following theorem suggests that the conditions `$\varphi$ is
monotone increasing and $\lim \limits_{t\rightarrow \infty}
\varphi(t)=\infty $' of Theorem 1.5 are redundant.

\smallskip

THEOREM 2.1. Let ~$T,~f:X\rightarrow X$ be weakly compatible
selfmaps. If $T$ is weakly contractive with respect to $f$ such that
$T(X)\subseteq f(X)$ and $f(X)$ is a complete subspace of $X$, then
$f$ and $T$ have a unique common fixed point in $X$.

\smallskip
\smallskip

PROOF. Let $x_{0}\in X.$ Since $T(X)\subseteq f(X)$,
 there exists $x_{1}\in X$ such that $Tx_{0}=fx_1.$
On continuing this process, inductively we get a sequence
$\{x_{n}\}$ in $X$ such that $y_{n}=T(x_{n})=f(x_{n+1})$.

 We now
show that the sequence $\{d(fx_n,fx_{n+1})\}$ is a decreasing
sequence.\\Now consider,
\begin{equation} d(fx_{n+1},fx_{n+2})=d(T(x_{n}),T(x_{n+1})) \leq
d(f(x_{n}),f(x_{n+1}))-\varphi (d(f(x_{n}),f(x_{n+1})). \label{2}
\end{equation}

Hence,
 \begin{equation} d(fx_{n+1},fx_{n+2})\leq d(f(x_{n}),f(x_{n+1}))~ \text{for all}~ n=0,1,2,\cdots. \label{3}
\end{equation}

Hence the sequence $\{d(f(x_{n}),f(x_{n+1}))\}$ is a decreasing
sequence of non-negative reals and converges to a limit $l$ (say)
and $l\geq 0$.

We claim that $l=0$. Suppose $l> 0$.

Letting $n\rightarrow \infty$ in (\ref{2}), by the continuity of
$\varphi$, we get $l\leq l-\varphi(l)$, a contradiction. Hence,
$l=0$. i.e., \begin{equation} \lim \limits_{n\rightarrow
\infty}d(f(x_{n}),f(x_{n+1}))=0. \label{4} \end{equation}

We now claim that $\{y_{n}\}$ is Cauchy. By (\ref{3}) and (\ref{4}),
it is sufficient to show that $\{y_{2n}\}$ is Cauchy. Otherwise,
there exists an $\varepsilon > 0$ and there exist sequences
$\{m_{k}\}$ and $\{n_{k}\}$ with $m_{k}> n_{k}> k$ such that
\begin{equation} d(y_{2m_{k}},y_{2n_{k}})\geq \varepsilon~~~~
\text{and}~~~~~ d(y_{2m_{k}-2},y_{2n_{k}})<\varepsilon. \label{5}
 \end{equation}
Hence, \begin{equation}\varepsilon \leq \liminf
\limits_{k\rightarrow \infty}d(y_{2m_{k}},y_{2n_{k}}). \label{6}
\end{equation}
For each positive integer $k$, by triangle inequality we have,
\[d(y_{2m_{k}},y_{2n_{k}})\leq d(y_{2m_{k}},y_{2m_{k}-1})+
d(y_{2m_{k}-1},y_{2m_{k}-2}) + d(y_{2m_{k}-2},y_{2n_{k}})\] On
taking limit supremum of both sides, as $k\rightarrow \infty$, we
get
\begin{equation}\limsup\limits_{k\rightarrow
\infty}d(y_{2m_{k}},y_{2n_{k}})\leq\varepsilon .\label{7}
\end{equation}
Hence, from (\ref{6}) and (\ref{7}), we have
\begin{equation}
\lim \limits_{k \to \infty}{d(y_{2m_{k}},y_{2n_{k}})}=\varepsilon.
\label{8}\end{equation} Now \[d(y_{2m_{k}},y_{2n_{k}+1})\leq
d(y_{2m_{k}},y_{2n_{k}})+ d(y_{2n_{k}},y_{2n_{k}+1})\] On taking
limit supremum, as $k\rightarrow \infty$, we get
\begin{equation}
\limsup\limits_{k\rightarrow
\infty}d(y_{2m_{k}},y_{2n_{k}+1})\leq\varepsilon. \label{9}
\end{equation}
Again  we have \[d(y_{2m_{k}},y_{2n_{k}})\leq
d(y_{2m_{k}},y_{2n_{k}+1})+ d(y_{2n_{k}+1},y_{2n_{k}}).\] On taking
limit infimum, as $k\rightarrow \infty$, we get
\begin{equation}
\varepsilon \leq\liminf\limits_{k\rightarrow
\infty}d(y_{2m_{k}},y_{2n_{k}+1}).\label{10} \end{equation}
 From  (\ref{9}) and (\ref{10}), we have
\begin{equation}
\lim \limits_{k \to \infty}{d(y_{2m_{k}},y_{2n_{k}+1})}=\varepsilon.
\label{11}
\end{equation}
Similarly, we can show that
\begin{equation}
\lim \limits_{k \to
\infty}{d(y_{2m_{k}-1},y_{2n_{k}})}=\varepsilon.\label{12}\end{equation}
Now consider,
\[{d(y_{2m_{k}},y_{2n_{k}+1})}= {d(Tx_{2m_{k}},Tx_{2n_{k}+1})}~~~~~~~~~~~~~~~~~~~~~~~\]
\[~~~~~~~~~~~~~~~~~~~~~~~~~~~~~\leq
{d(fx_{2m_{k}},fx_{2n_{k}+1})}-\varphi({d(fx_{2m_{k}},fx_{2n_{k}+1})})\]
\[ ~~~~~~~~~~~~~~~~~~~~~~~~~~~~~~~\leq
{d(Tx_{2m_{k}-1},Tx_{2n_{k}})}-\varphi({d(Tx_{2m_{k}-1},Tx_{2n_{k}})}).\]
\begin{equation}
~~~~~~~~~~~~~~~~~~~~~~\leq
{d(y_{2m_{k}-1},y_{2n_{k}})}-\varphi({d(y_{2m_{k}-1},y_{2n_{k}})}).
\label{13}
\end{equation}
Letting $k\rightarrow \infty$ in (\ref{13}), using (\ref{11}),
(\ref{12}) and the continuity of $\varphi$, we get $\varepsilon\leq
\varepsilon - \varphi(\varepsilon)$, a contradiction. Hence,
$\{y_{2n}\}$ is Cauchy so that $\{y_{n}\}$ is a Cauchy sequence in
$X$. Thus $\{fx_{n+1}\}$ is a Cauchy sequence in $X$. Since $f(X)$
is complete and $\{fx_{n+1}\}\subset f(X)$, we have
\begin{equation}
\lim\limits_{n\rightarrow \infty}fx_{n+1}=fu,~ \text{for some } u\in
X.\label{14}
\end{equation}
Next we claim that $T(u)=f(u)$.
Now consider,
\begin{equation}
d(fx_{n+1},Tu)=d(Tx_{n},Tu)\leq
  d(fx_{n},fu)-\varphi(d(fx_{n},fu)). \label{15}\end{equation}
Letting $n\rightarrow \infty$, from (\ref{14}) using (\ref{15}) and
the continuity of $\varphi$, we get \[d(fu,Tu)\leq
  d(fu,fu)-\varphi(d(fu,fu))=0.\]

Hence, $Tu=fu=z$ (say). Since the pair of maps $(f,T)$ is weakly
compatible, we have $Tfu=fTu$ and hence $Tz=fz$.

We now claim that $Tz=z$. Suppose $Tz\neq z$. Consider
\[d(Tz,z)=d(Tz,Tu)\leq d(fz,fu)-\varphi(d(fz,fu))\]
\[~~~~~~~~~~~~~~~~~~~~~~~~
=d(Tz,z)-\varphi(d(Tz,z)),\] a contradiction. Hence, $Tz=z$.

Hence, $fz=Tz=z$. The uniqueness of $z$ follows from the weakly
contractive nature of $T$.
 Hence, the theorem follows.

\smallskip

The following is an example in support of our Theorem 2.1.
\smallskip

EXAMPLE 2.2. Let $X=R_{+}$ with the usual metric. We define mappings
$f$ and $T$ on $X$ by
\[\begin{array}{lr}
 f(x)=\left\{\begin{array}{ll}
                  \frac{1}{3} & \text{if}~0\leq x < \frac{2}{3} \\\\
                 \frac{2}{3} & \text{if}~x=\frac{2}{3} \\\\
                 \frac{5}{6} & \text{if}~x>\frac{2}{3}\end{array}\right.\nonumber &\mbox{}\text{and} ~~~~  T(x)=\left\{\begin{array}{ll}
                  \frac{5}{6} & \text{if}~0\leq x < \frac{2}{3} \\\\
                 \frac{2}{3} & \text{if} ~x\geq \frac{2}{3}.\end{array}\right.\nonumber
\end{array}\]

\smallskip
 Clearly, the pairs $(f,T)$ is weakly compatible, $T(X)\subseteq f(X)$ and $f(X)=\{\frac{1}{3}, \frac{2}{3},\frac{5}{6}\}$ is complete.

We define $\varphi:R_{+}\rightarrow R_{+}$ by
$\varphi(t)=\left\{\begin{array}{ll}
              \frac{3}{2}t^{2} & \text{if}~0\leq t\leq \frac{1}{3} \\
              \frac{2}{9(1+t)} & \text{if}~t\geq\frac{1}{3}.
            \end{array}\right.\nonumber$ Clearly $\varphi\in \Phi$. With this $\varphi$, $T$ is weakly contractive with respect to $f$.

Hence, $f$ and $T$ satisfy all the conditions of Theorem
 2.1 and $\frac{2}{3}$ is the unique common fixed point of $f$ and
 $T$.

Here  $\varphi$ is not monotonically increasing on $R_{+}$, and
$\lim \limits_{n\rightarrow \infty}\varphi(t)\neq \infty$  so that
Theorem 1.5 is not applicable.
\section{Main Result}

\smallskip
\hskip 0.6cm THEOREM 3.1. Let $(X,d)$ be a metric space and let
$T,~f:X\rightarrow X$ be weakly compatible selfmaps satisfying
property (E. A). Assume that $T$ is weakly contractive with respect
to $f$. If $f(X)$ is closed, then $f$ and $T$ have a unique common
fixed point in $X$.

\smallskip
\smallskip

PROOF. Since the pair $(f,T)$ satisfies property (E. A), there
exists a sequence $\{x_{n}\}$ in $X$ such that$\lim\limits_{n\to
\infty}fx_{n}=\lim\limits_{n\to \infty}Tx_{n}=z,$ for some $z$ in
$X$.

Since $f(X)$ is closed, $z=f(u)$ for some $u\in X$. \\Now replacing
$u$ for $x$ and $x_{n}$ for $y$ in (\ref{1}), we get
\begin{equation}
d(Tu,Tx_{n})\leq d(fu,fx_{n})-\varphi(d(fu,fx_{n})). \label{3.1}
\end{equation}
Letting $n\rightarrow \infty$ in (\ref{3.1}) , by the continuity of
$\varphi$, we get \[d(Tu,z)\leq d(fu,z)-\varphi(d(fu,z))=0.\] Hence,
$Tu=z$.\\Hence,
\begin{equation}
fu=Tu=z. \label{3.2}
\end{equation}

Since $f$ and $T$ are weakly compatible, from (\ref{3.2}), we have
$fz=Tz$. \\If $z\neq Tz$, then from the inequality (\ref{1}), we
have \[d(Tz,z)=d(Tz,Tu)\leq d(fz,fu)-\varphi(d(fz,fu))\]
\[~~~~~~~~~~~~~~~~~~~~~~~~~~~~=  d(Tz,Tu)-\varphi(d(Tz,Tu)),\]
a contradiction. Hence, $Tz=z$.

Hence, $fz=Tz=z$.

The uniqueness of $z$ follows from the inequality (\ref{1}).

This complete the proof the Theorem.

\smallskip

Since two noncompatible selfmaps of a metric space $(X,d)$ satisfy
the property \\(E. A), we get the following corollary.

\smallskip

COROLLARY 3.2. Let $(X,d)$ be a metric space and let
$f,T:X\rightarrow X$ be noncompatible and weakly compatible
selfmaps. Assume that $T$ is weakly contractive with respect to $f$.
If $f(X)$ is closed, then $f$ and $T$ have a unique common fixed
point in $X$.

\smallskip

In Theorem 3.1, if we choose $\varphi \in \Phi$ such that
$\varphi(t)=t-\psi(t)$, where $\psi \in \Phi$ with $\psi(t)<t$ for
$t>0$, we get the following Corollary.
\smallskip

COROLLARY 3.3. Let $(X,d)$ be a metric space and let
$T,~f:X\rightarrow X$ be weakly compatible selfmaps satisfying
property (E. A). Assume that there exists a $\psi \in \Phi$ with
$\psi(t)<t$ for $t>0$ such that
\begin{equation}
d(Tx,Ty)\leq \psi(d(fx,fy))~ \text{for all} ~x,y \in X. \label{3.3}
\end{equation}
If $f(X)$ is closed, then $f$ and $T$ have
a unique common fixed point in $X$.

\smallskip

In this case, when $T$ and $f$ satisfy the inequality (\ref{3.3}),
we say that \emph{$T$ is a Boyd-Wong type contraction with respect
to $f$} \cite{Boyd}.

The following is an example in support of Theorem 3.1.
\smallskip

EXAMPLE 3.4. Let $X=[\frac{1}{2},1]$ with the usual metric. We
define mappings $f$ and $T$ on $X$ by
\[\begin{array}{lr}
 f(x)=\left\{\begin{array}{ll}
                  1 & \text{if}~\frac{1}{2}\leq x < \frac{2}{3} \\\\
                 x & \text{if} ~\frac{2}{3}\leq
  x\leq 1 \end{array}\right.\nonumber &\mbox{}\text{and} ~~~~  T(x)=\left\{\begin{array}{ll}
                  ~~\frac{1}{2} & \text{if}~\frac{1}{2}\leq x < \frac{2}{3} \\\\
                 1-\frac{1}{2}x & \text{if }~\frac{2}{3}\leq
  x\leq 1 .\end{array}\right.\nonumber
\end{array}\]

\smallskip

Since the sequence $\{x_{n}\}$, $x_{n}=\frac{2}{3}+\frac{1}{n},
~n\geq 4$, in $X$ with $\lim\limits_{n\to
\infty}fx_{n}=\lim\limits_{n\to
 \infty}Tx_{n}=\frac{2}{3}$,  the pair $(f,T)$
 satisfy property (E. A).

 Clearly, the pairs $(f,T)$ is weakly compatible and $f(X)=[\frac{2}{3}, 1]$ is closed.

We define $\varphi:R_{+}\rightarrow R_{+}$ by
$\varphi(t)=\left\{\begin{array}{ll}
              \frac{3}{2}t^{2} & \text{if}~0\leq t\leq \frac{2}{3} \\
              \frac{10}{9(1+t)} & \text{if}~t\geq\frac{2}{3}.
            \end{array}\right.\nonumber$ Clearly $\varphi\in \Phi$. With this $\varphi$, $T$ is weakly contractive with respect to $f$.

Hence, $f$ and $T$ satisfy all conditions of Theorem
 3.1 and $\frac{2}{3}$ is the unique common fixed point of $f$ and
 $T$.

Further we mention that the pair $(f,T)$ is not compatible, for
$\lim\limits_{n\to
 \infty}d(fTx_{n},Tfx_{n})=\frac{1}{6}\neq 0$.

Here we observe that neither $T(X)\subseteq f(X)$ nor $f(X)\subseteq
T(X)$, and $\varphi$ is not monotonically increasing on $R_{+}$, so
that Theorem 1.5 and Theorem 2.1 are not applicable.

\smallskip

EXAMPLE 3.5. Let $l_{\infty}$ be the set of all bounded nonnegative
real numerical sequences $\{x_{n}\}$. We define metric $d$ on
$l_{\infty}$ by $d(x,y)=\sup \{|x_{n}-y_{n}|:n\in N \}$, where
$x=\{x_{n}\}$ and $y=\{y_{n}\}$ in $l_{\infty}$. Then
$(l_{\infty},d)$ is a complete metric space.

We define $T:l_{\infty}\rightarrow l_{\infty}$ by
$T(\{x_{n}\})=\{\frac{x_{n}}{1+x_{n}}\}$ and $f=I$, the identity map
on $l_{\infty}$.

\smallskip

 Let $x=\{x_{n}\}, y=\{y_{n}\} \in l_{\infty}$. Then
\[ d(Tx,Ty)=d(\{\frac{x_{n}}{1+x_{n}}\},\{\frac{y_{n}}{1+y_{n}}\})\]
\[ ~~~~~~~~~~~~~~~~~~~~~~~~~=\sup
\{|\frac{x_{n}}{1+x_{n}}-\frac{y_{n}}{1+y_{n}}|:n\in N\}\]
\[~~~~~~~~~~~~~~~~~~~\leq \sup\{\frac{|x_{n}-y_{n}|}{1+|x_{n}-y_{n}|}:n\in
N\}\]
 \[ ~~~~~~~~~~~~~~~~~~~\leq \frac{\sup \{|x_{n}-y_{n}|:n\in N\}}{1+\sup \{|x_{n}-y_{n}|:n\in
 N\}}\]
\[ ~~~~~~~~~~~~~~~~~~~~~~~~~~= d(\{x_{n}\},\{y_{n}\})-\varphi(
d(\{x_{n}\},\{y_{n}\}))\]
 \[ ~~~~~~~~~= d(x,y)-\varphi( d(x,y)),\] where $\varphi(t)=\frac{t^{2}}{1+t},~ t\geq
 0.$

Thus $T$ is a weakly contractive map with respect to $f$.

Clearly, $f$ and $T$ satisfy property (E.A), by choosing the
sequence $\{x_{n}\}, x_{n}=(0,0,\cdots)\in l_{\infty}$ for all
$n=1,2,\cdots$. Since the null vector is the only coincidence point
of $f$ and $T$ , we have $f$ and $T$ are weakly compatible on
$l_{\infty}$. Further, $f(l_{\infty})=l_{\infty}$ is closed. Hence,
$f$ and $T$ satisfy all conditions of Theorem 3.1 and the null
vector is the unique common fixed point of $f$ and $T$ in
$l_{\infty}$.

We observe that $T$ is not a contraction on $l_{\infty}$, for taking
$x=(0,0,\cdots)$,  for all $k\in(0,1)$ there exists
$y=(t,0,0,\cdots) \in l_{\infty}$, where $t\in (0,\frac{1-k}{k})$,
such that $d(Tx,Ty)=\frac{t}{1+t}> k~t=k~d(x,y)$.
\smallskip

COROLLARY 3.6. Let $K$ be a nonempty closed subset of a metric space
$(X,d)$ and $T:X\rightarrow X$. Assume that there exists a $\varphi
\in \Phi$ such that \[d(Tx,Ty)\leq d(x,y)-\varphi(d(x,y)),\]  for
all $x,y\in X$. If there exists a sequence $\{x_{n}\}$ in $K$ such
that $\lim \limits_{n\rightarrow \infty}Tx_{n}=\lim
\limits_{n\rightarrow \infty}x_{n}=z,~ z\in K$, then $z$ is the
fixed point of $T$ in $K$.

\smallskip

The following is an example in support of Corollary 3.6.

\smallskip

EXAMPLE 3.7. Let $X=R_{+}$ with the usual metric and $K=\{0\}\cup
[\frac{1}{3},1]$. We define a mapping $T$ on $K$ by
$Tx=\left\{\begin{array}{ll}
              ~~\frac{1}{3} & \text{if}~~~ x=0 \\
              \frac{1}{3}+\frac{1}{2}x & \text{if}~\frac{1}{3}\leq x\leq 1.
            \end{array}\right.\nonumber$

\smallskip
Then $T$ satisfies all the conditions of Corollary 3.6 with
 $\varphi:R_{+}\rightarrow R_{+}$ defined by
 $\varphi(t)=\frac{t^{2}}{1+t},~t\geq 0$ and $\frac{2}{3}$ is the
 unique fixed point of $T$.

We observe that the sequence $\{x_{n}\}$,
$x_{n}=\frac{2}{3}+\frac{1}{n},~n\geq 4$, is in $K$ with $\lim
\limits_{n\rightarrow \infty}Tx_{n}=\lim \limits_{n\rightarrow
\infty}x_{n}=\frac{2}{3}$.
\smallskip

REMARK 3.8. If we delete the condition `$f(X)$ is closed' from
Theorem 2.1, then the maps $f$ and $T$ may have no common fixed
points, which is shown by the following example.

\smallskip
\smallskip

EXAMPLE 3.9. Let $X=[\frac{1}{2},1]$ with the usual metric. We
define mappings $f$ and $T$ on $X$ by

\[\begin{array}{lr}
 f(x)=\left\{\begin{array}{ll}
                  1 & \text{if}~\frac{1}{2}\leq x \leq\frac{2}{3} \\\\
                 x & \text{if} ~\frac{2}{3}<
  x\leq 1 \end{array}\right.\nonumber &\mbox{}\text{and} ~~~~  T(x)=\left\{\begin{array}{ll}
                  ~~~\frac{1}{2} & \text{if}~\frac{1}{2}\leq x \leq \frac{2}{3} \\\\
                 1-\frac{1}{2}x & \text{if} ~\frac{2}{3}<
  x\leq 1. \end{array}\right.\nonumber
\end{array}\]

\smallskip

Since the sequence $\{x_{n}\}$, $x_{n}=\frac{2}{3}+\frac{1}{n},~
n\geq 4$,
 in $X$ with $\lim\limits_{n\to \infty}fx_{n}=\lim\limits_{n\to
 \infty}Tx_{n}=\frac{2}{3}$,  the pair $(f,T)$
 satisfy property (E. A). Clearly, the pairs $(f,T)$ is weakly compatible.

 We define
 $\varphi:R_{+}\rightarrow R_{+}$ by
 $\varphi(t)=\frac{3}{2}t^{2},~t\geq 0$. Clearly $\varphi \in \Phi$.
 With this $\varphi$, $T$ is weakly contractive with respect to $f$.

But  $f(X)=(\frac{2}{3}, 1]$ is not closed. We observe that $f$ and
$T$ have no common fixed points.
\smallskip

REMARK 3.10. In Theorem 3.1, if we relax the condition `$f$ and $T$
satisfy property (E. A)' then they may have no common fixed points.

\smallskip
\smallskip

EXAMPLE 3.11. Let $X=[\frac{1}{2},1]$ with the usual metric. We
define mappings $f$ and $T$ on $X$ by

\[\begin{array}{lr}
 f(x)=\left\{\begin{array}{ll}
                  1 & \text{if}~\frac{1}{2}\leq x \leq\frac{2}{3} \\\\
                 \frac{2}{3} & \text{if} ~\frac{2}{3}<
  x\leq 1 \end{array}\right.\nonumber &\mbox{}\text{and} ~~~~  T(x)=\left\{\begin{array}{ll}
                  \frac{2}{3} & \text{if}~\frac{1}{2}\leq x \leq \frac{2}{3} \\\\
                 \frac{1}{2} & \text{if} ~\frac{2}{3}<
  x\leq 1. \end{array}\right.\nonumber
\end{array}\]

\smallskip

Here $f$ and $T$ are trivially weakly compatible on $X$ and
$f(X)=\{\frac{2}{3},1\}$ is closed. Further, $T$ is weakly
contractive with respect to $f$ with $\varphi(t)=\frac{1}{2}t,~t\geq
0$.

But $f$ and $T$ do not satisfy property (E. A), since for any
 sequence $\{x_{n}\}$ in $X$ we have $\lim\limits_{n\to \infty}fx_{n}\neq \lim\limits_{n\to
 \infty}Tx_{n}$ in $X$. We observe that $f$ and $T$ have no common fixed points.

REMARK 3.12. In Theorem 1.5, the authors assumed the condition
$T(X)\subseteq f(X)$, where as in the results of this paper, this
condition is relaxed by imposing the condition property (E. A).

\smallskip

\textbf{Acknowledgment.} The authors sincerely thank Prof. K. P. R.
Sastry for his valuable suggestions in the construction of Example
3.5.


\smallskip

\begin{thebibliography}{9}
\bibitem {Aamri} M. Aamri and D. El Moutawakil, Some new common fixed
point theorems under strict contractive conditions, J. Math. Anal.
Appl., 270(2002), 181-- 188.
\bibitem {Aliouche} A. Aliouche, Common fixed point theorems of
Greg$\ddot{u}$s type for weakly compatible mappings satisfying
generalized contractive conditions, J. Math. Anal. Appl., 341(2008),
707--719.
\bibitem {Beg} I. Beg and M. Abbas, Coincidence points and invariant approximation for
mappings satisfying generalized weak contractive condition, Fixed
point theory and applications, 2008, ID.74503, 1--7.
\bibitem {Boyd} D. W. Boyd and T. S. W. Wong, On nonlinear contractions, Proc. Amer. Math. Soc., 20 (1969), \\458 -- 464.
\bibitem {Imdad} M. Imdad and J. Ali, Jungck's common fixed point
theorem and E. A property, Acta Mathematica Sinica, English series,
Jan., 2008, Vol. 24, No.1, pp. 87--94

\bibitem{Jungck1} G. Jungck,
Compatible mappings and common fixed points,  Internat. J. Math.
Math. Sci.,  9 (4) (1986), 771--779.
\bibitem{Jungck2}  G. Jungck and B.E. Rhoades, Fixed point for
set-valued functions without continuity, Indian J. Pure and Appl.
Math. 29(3)( 1998), 227--238.
\bibitem {Liu} W. Liu, J. Wu and Z. Li, Common fixed points of single-valued and multi-valued
maps, Int. J. Math. Math. Sci., 19(2005), 3045--3055.

\bibitem {Pathak} H.K. Pathak, Rozana Rodriguez-Lopez and R.K.
Verma, A common fixed point theorem using implicit relation and
property (E.A) in metric spaces, FILOMAT., 21(2)(2007), 211--234.

\end{thebibliography}
\end{document}